\documentclass{amsart}

\newcommand{\Comp}{\mathcal{C}omp}
\newcommand{\e}{\varepsilon}
\newcommand{\IR}{\mathbb R}
\newcommand{\Reg}{\mathrm{Reg}}

\newcommand{\supp}{\mathrm{supp}}

\newcommand{\w}{\omega}
\newcommand{\IN}{\mathbb N}
\newcommand{\pr}{\mathrm{pr}}

\newtheorem{theorem}{Theorem}[section]
\newtheorem{proposition}[theorem]{Proposition}
\newtheorem{corollary}[theorem]{Corollary}
\newtheorem{problem}[theorem]{Problem}

\begin{document}
\title[Clifford subsemigroups of functor-semigroups]{Characterizing compact Clifford semigroups that embed into convolution and functor-semigroups}
\author[T.~Banakh, M.~Cencelj, O.~Hryniv, and D.~Repov\v s]{Taras Banakh, Matija Cencelj, Olena Hryniv, and Du\v san Repov\v s}
\begin{abstract}
 We study algebraic and topological properties of the
convolution semigroups of probability measures on
a topological groups and show that a compact Clifford
topological semigroup $S$ embeds into the convolution semigroup
$P(G)$ over some topological group $G$ if and only if  $S$ embeds
into the semigroup $\exp(G)$ of compact subsets of $G$ if and only if $S$ is an inverse semigroup and has zero-dimensional maximal semilattice. We
also show that such a Clifford semigroup $S$ embeds into the
functor-semigroup $F(G)$ over a suitable compact topological group $G$ for each weakly normal monadic functor $F$ in the
category of compacta such that $F(G)$ contains a $G$-invariant
element (which is an analogue of the Haar measure on $G$).
\end{abstract}

\address{Instytut Matematyki, Akademia \'Swi\c etokrzyska, Kielce, Poland \newline and
 Department of Mathematics, Ivan Franko National University of Lviv, Lviv, Ukraine}
\email{tbanakh@yahoo.com}

\address{Institute of Mathematics, Physics
and Mechanics, and Faculty of Education, University
of Ljubljana, P.O.B. 2964, Ljubljana, 1001, Slovenia} \email
{matija.cencelj@guest.arnes.si}

\address{Department of Mathematics, Ivan Franko National University of Lviv, Lviv, Ukraine}
\email{olena\_hryniv@ukr.net}

 \address{Faculty of Mathematics and Physics,
 and Faculty of Education,
 University
of Ljubljana,  P.O.B. 2964, Ljubljana, 1001, Slovenia}
\email
{dusan.repovs@guest.arnes.si}

\date{\today}

 \maketitle
\section{Introduction}
According to \cite{Ber} (and \cite{Trn}) each (commutative)
semigroup $S$ embeds into the global semigroup $\Gamma(G)$ over a
suitable (abelian) group $G$. The global semigroup $\Gamma(G)$
over $G$ is the set of all non-empty subsets of $G$ endowed with
the semigroup operation $(A,B)\mapsto AB=\{ab:a\in A,\; b\in B\}$.
If $G$ is a topological group, then the global semigroup
$\Gamma(G)$ contains a subsemigroup $\exp(G)$ consisting of all
non-empty compact subsets of $G$ and carrying a natural topology
which makes it a topological semigroup. This is the Vietoris
topology generated by the sub-base consisting of the sets
$$U^+=\{K\in\exp(G):K\subset U\}\mbox{ and
}U^-=\{K\in\exp(G):K\cap U\ne\emptyset\}$$ where $U$ runs over
open subsets of $G$. Endowed with the Vietoris topology the
semigroup $\exp(G)$ will be referred to as the {\em
hypersemigroup} over $G$ (because its underlying topological space
is the hyperspace $\exp(G)$ of $G$, see \cite{TZ}). The problem of
detecting topological semigroups embeddable into the
hypersemigroups over topological groups has been considered in the
literature, see \cite{Ber}.

 This problem was resolved in \cite{BHr} for the class of Clifford compact topological semigroups: such a semigroup $S$ embeds into the hypersemigroup over a topological group if and only if
the set $E$ of idempotents of $S$ is a zero-dimensional commutative
subsemigroup of $S$. This
characterization implies the result of  \cite{BL} that  the closed interval $[0,1]$ with the operation of the minimum does
not embed into the hypersemigroup over a topological group.

We recall that a semigroup $S$ is {\em Clifford} if $S$ is the union of
its subgroups. We say that a topological semigroup $S_1$ embeds
into another topological semigroup $S_2$ if there is a semigroup
homomorphism $h:S_1\to S_2$ which is a topological embedding.

In this paper we shall apply the already mentioned result of \cite{BHr}
and  shall characterize Clifford compact semigroups embeddable into
the convolution semigroups $P(G)$ over topological groups $G$. The
convolution semigroup $P(G)$ consists of probability Radon
measures on $G$ and carries the $*$-weak topology generated by the
sub-base $\{\mu\in P(G):\mu(U)>a\}$ where $a\in\IR$ and $U$ runs
over open subsets of $G$. A measure $\mu$ defined on the
$\sigma$-algebra of Borel subsets of $G$ is called {\em Radon} if
for every $\e>0$ there is a compact subset $K\subset G$ with
$\mu(K)>1-\e$. The semigroup operation on $P(G)$ is given by the
convolution measures. We recall that the {\em convolution}
$\mu*\nu$ of two measures $\mu,\nu$ is the measure assigning to
each bounded continuous function $f:G\to\IR$ the value of the integral
$\int_{\mu*\nu}f=\int_\eta\int_\mu f(xy)dydx$. For more detail information on the convolution semigroups, see \cite{Heyer}, \cite{Par}.

The following theorem is the principal result of this paper.

\begin{theorem}\label{t1} For any Clifford compact topological semigroup $S$ the following assertions are equivalent:
\begin{enumerate}
\item $S$ embeds into the hypersemigroup $\exp(G)$ over a topological group $G$;
\item $S$ embeds into the convolution semigroup $P(G)$ over a topological group $G$;
\item The set $E$ of idempotents of $S$ is a zero-dimensional commutative subsemigroup of $S$.
\end{enumerate}
\end{theorem}

This theorem will be applied to a characterization of Clifford compact
topological semigroups embeddable into the hyperpsemigroups or
convolution semigroups over topological groups $G$ belonging to
certain varieties of topological groups. A class $\mathcal G$ of
topological groups is called a {\em variety} if it is closed under
 arbitrary Tychonov products, and taking closed subgroups, and
quotient groups by closed normal subgroups.

\begin{theorem}\label{t2} Let $\mathcal G$ be a non-trivial variety of topological groups. For a Clifford compact topological semigroup $S$
the following assertions are equivalent:
\begin{enumerate}
\item $S$ embeds into the hypersemigroup $\exp(G)$ over a topological group $G\in\mathcal G$;
\item $S$ embeds into the convolution semigroup $P(G)$ over a topological group $G\in\mathcal G$;
\item The set $E$ of idempotents is a zero-dimensional commutative subsemigroup of $S$ and all
closed subgroups of $S$ belong to the class $\mathcal G$.
\end{enumerate}
\end{theorem}

In fact, the equivalence of the first and last statements in
Theorems~\ref{t1} and \ref{t2} was proved in Theorems 3 and 4 of \cite{BHr} so
it remains to prove the equivalence of the assertions (1) and (2). This
will be done in Proposition~\ref{pr1}. 

We recall that a semigroup $S$ is
called {\em regular} if each element $x\in S$ is {\em regular} in
the sense that $xyx=x$ for some $y\in S$. An element $x\in S$ is
called ({\em uniquely}) {\em invertible} if there is a (unique)
element $x^{-1}\in S$ (called the {\em inverse} of $x$) such that
$xx^{-1}x=x$ and $x^{-1}xx^{-1}=x^{-1}$. A semigroup $S$ is called
{\em inverse} if each element of $S$ is uniquely invertible. By
\cite[1.17]{CP}, \cite[II.1.2]{Pet} a semigroup $S$ is inverse if
and only if it is regular and the set $E$ of idempotents of $S$ is
a commutative subsemigroup of $S$. An inverse semigroup $S$ is
Clifford if and only if $xx^{-1}=x^{-1}x$ for all $x\in S$. In
this case $S$ decomposes into the union $S=\bigcup_{e\in E}H_e$ of
the maximal subgroups $H_e=\{x\in S:xx^{-1}=e=x^{-1}x\}$ of $S$
parametrized by idempotents $e$ of $S$.

We recall that a topological semigroup $S$ is called a {\em
topological inverse semigroup} if $S$ is an inverse semigroup and
the inversion map $(\cdot)^{-1}:S\to S$, $(\cdot)^{-1}:x\mapsto
x^{-1}$ is continuous. The set $E$ of idempotents of a topological
inverse semigroup $S$ is a closed commutative subsemigroup of $S$
called the {\em idempotent semilattice} of $S$. We say that two
idempotents $e,f\in E$ are {\em incomparable} if their product $ef$
differs from $e$ and $f$. Two elements $x,y$ of an inverse
semigroup $S$ are called {\em conjugate} if $x=zyz^{-1}$ and
$y=z^{-1}xz$ for some element $z\in S$. For any idempotent $e\in E$
let ${\uparrow}e=\{f\in E:ef=e\}$ denote the principal filter of
$e$. A topological space $X$ is called {\em totally disconnected}
if for any distinct points $x,y\in X$ there is a closed-and-open
subset $U\subset X$ containing $x$ but not $y$.

The following proposition shows that the semigroups $\exp(G)$ and
$P(G)$ over a topological group $G$ have the same regular
subsemigroups (which are necessarily  topological inverse
semigroups). Moreover, regular subsemigroups of $\exp(G)$ or
$P(G)$ have many specific topological and algebraic features.

\begin{proposition}\label{pr1} Let $G$ be a topological group. A topological regular semigroup $S$ embeds into $P(G)$
if and only if $S$ embeds into $\exp(G)$. If the latter happens, then
\begin{enumerate}
\item $S$ is a topological inverse semigroup;
\item The idempotent semilattice $E$ of $S$ has totally disconnected principal filters ${\uparrow}e$, $e\in E$;
\item An element $x\in S$ is an idempotent if and only if $x^2x^{-1}$ is an idempotent;
\item Any distinct conjugated idempotents of $S$ are incomparable.
\end{enumerate}
\end{proposition}

This proposition allows one to construct many examples of
topological regular semigroups non-embeddable into the
hypersemigroups or convolution semigroups over a topological
groups. The first two assertions of this proposition imply the result
of \cite{BL} to the effect that non-trivial semigroups of left (or right) zeros as well as connected
topological semilattices do not embed into the hypersemigroup
$\exp(G)$ over a topological group $G$. The last two assertions imply
that the semigroups $\exp(G)$ and $P(G)$ do not contain Brandt
semigroups and bicyclic semigroups. By a {\em Brandt semigroup} we
understand a semigroup of the form $B(H,I)=I\times H\times
I\cup\{0\}$ where $H$ is a group, $I$ is a non-empty set, and the
product $(\alpha,h,\beta)*(\alpha',h',\beta')$ of two non-zero
elements of $B(H,I)$ is equal to $(\alpha,hh',\beta')$ if
$\beta=\alpha'$ and $0$ otherwise. A {\em bicyclic semigroup} is a
semigroup generated by two elements $p,q$ with the relation
$qp=1$. Brandt semigroups and byciclic semigroups play an important role in
the structure theory of inverse semigroups, see \cite{Pet}.

In fact, the semigroups $\exp(G)$ and $P(G)$ are special cases
of the so-called functor-semigroups introduced by Teleiko and
Zarichnyi \cite{TZ}. They observed that any weakly normal monadic
functor $F:\Comp\to\Comp$ in the category of compact Hausdorff
spaces lifts to the category of compact topological semigroups,
which means that for any compact topological semigroup $X$ the
space $FX$ possesses a natural semigroup structure. The semigroup
operation $*$ on $FX$ can be defined by the following formula
$$a* b=Fp(a\otimes b)\mbox{ for }a,b\in FX$$
where $p:X\times X\to X$ is the semigroup operation of $X$ and
$a\otimes b\in F(X\times X)$ is the tensor product of the elements
$a,b\in FX$, see \cite[\S3.4]{TZ}.

Therefore we actually consider in this paper the following general problem:

\begin{problem} Given a weakly normal monadic functor $F:\Comp\to\Comp$,  find a characterization of compact (regular, inverse, Clifford)
topological semigroups embeddable into the semigroup $FX$ over
a compact topological group $X$. Given a compact topological group $X$ describe
invertible elements and idempotents of the semigroup $FX$.
\end{problem}

Observe that for the functors $\exp$ and $P$ the answer to the
first part of this problem is given in Theorem~\ref{t1}.
Functor-semigroups induced by the functors $G$ of inclusion hyperspaces and $\lambda$ of superextension have been studied in  \cite{BGN}--\cite{BG4},\cite{G1}.

 In fact, Theorem~\ref{t2}  also can be partly generalized to some
monadic functors $F$ (including the functors $\exp$, $P$, $G$ and $\lambda$). Given a compact topological group $G$ let us
define an element $a\in F(G)$ to be {\em $G$-invariant} if
$g*a=a=a*g$ for every $g\in G$. Here we identify $G$ with a
subspace of $F(G)$ (which is possible because $F$, being weakly
normal, preserves singletons). A $G$-invariant element in $F(G)$ exists for the functors $\exp$, $P$, $\lambda$, and $G$. For the functors $\exp$ and $P$ a $G$-invariant element on $F(G)$ is unique: it is $G\in\exp(G)$ and the Haar measure on $G$, respectively.

\begin{theorem}\label{t3} Let $F:\Comp\to\Comp$ be a  weakly normal monadic functor such that for every compact topological group $G$
the semigroup $F(G)$ contains a $G$-invariant element. Each Clifford compact topological inverse semigroup $S$ with zero-dimensional
idempotent semilattice $E$ embeds into the functor-semigroup $F(G)$ over the compact topological group $G=\prod_{e\in E}\widetilde H_e$
where each $\widetilde H_e$ is a non-trivial compact topological group containing the maximal subgroup $H_e\subset S$ corresponding to an idempotent $e\in E$ of $S$.
\end{theorem}

\begin{proof}
By Theorem 3 of \cite{BHr} (see also \cite{Hr}), each Clifford compact topological
inverse semigroup $S$ with zero-dimensional idempotent semilattice
$E$ embeds into the product $\prod_{e \in E}H_e ^0$, where $H_e ^0$ stands for the extension of the maximal subgroup $H_e$ by an isolated point
$0 \notin H_e$ such that $x0=0x=0$ for all $x\in H_e$. For every
idempotent $e\in E$, fix a non-trivial compact topological
group $\widetilde{H}_e$ containing $H_e$. By our hypothesis, the space $F(\widetilde{H}_e)$ contains an $\widetilde{H}_e$-invariant element $z_e\in F(\widetilde{H}_e)$. Then $H_e^0$ can be identified with the closed subsemigroup $H_e\cup\{z_e\}$ of $F(\widetilde{H}_e)$ and the product $\prod_{e\in E}H_e^0$ can be identified with a subsemigroup of the product $\prod_{e\in E}F(\widetilde H_e)$.
By \cite[p.126]{TZ}, the latter product can be identified with a subspace (actually a subsemigroup) of $F(\prod_{e\in E}\widetilde H_e)=F(G)$, where $G=\prod_{e\in E}\widetilde H_e$. In this way, we obtain an embedding of $S$ into $F(G)$.
\end{proof}

As we have said, the functors $\lambda$ of superextension and $G$ of inclusion hyperspaces satisfy the hypothesis of Theorem~\ref{t3}.
However, Proposition~\ref{pr1} is specific for the functor $P$ and cannot be generalized to the functors $\lambda$ or $G$.

Indeed, for the 4-element cyclic group $C_4$ the semigroup $\lambda(C_4)$ is isomorphic to the commutative inverse semigroup $C_4\oplus C_2^1$, where $C_2^1=C_2\cup\{1\}$ is the result of attaching an external unit to the 2-element cyclic group $C_2$, (see \cite{BGN}). On the other hand, the 12-element  semigroup $C_4\oplus C_2^1$ cannot be embedded into $\exp(C_4)$ because the set of regular elements of $\exp(C_4)$ consists of 7 elements (which are shifted subgroups of $C_4$).
Also the commutative inverse semigroup $\lambda(C_4)\cong C_4\oplus C_2^1$ can be embedded into $G(C_4)$ (because $\lambda$ is a submonad of $G$) but cannot embed into $\exp(C_4)$.

\section{Idempotents and invertible elements of the convolution semigroups}

In this section we prove Proposition~\ref{pr1}. For each
topological group $G$ the semigroups $P(G)$ and $\exp(G)$ are
related via the map of the support. We recall that the {\em support}
of a Radon measure $\mu\in P(G)$ is the closed subset
$$S_\mu=\{x\in G:\mu(Ox)>0\mbox{ for each neighborhood $Ox$ of $x$}\}$$of $G$.
Let $2^G$ denote the semigroup of all non-empty closed subsets of
$G$ endowed with the semigroup operation $A*B=\overline{AB}$. By
$$\supp:P(G)\to 2^G,\; \supp:\mu\mapsto S_\mu$$
we denote the support map.

The following proposition is well-known, see (the proof of) Theorem 1.2.1 in \cite{Heyer}.

\begin{proposition} Let $G$ be a topological group. For any measures $\mu,\nu\in P(G)$ the following holds:
$S_{\mu*\nu}=\overline{S_\mu\cdot S_\nu}$. This means that the support map $\supp:P(G)\to 2^G$ is a semigroup
homomorphism.
\end{proposition}

We shall show that for any regular element $\mu$ of the convolution semigroup $P(G)$ the support
$S_\mu$ is compact and thus belongs to the subsemigroup $\exp(G)$ of $2^G$.
First, we characterize idempotent measures on a topological group $G$.

A measure $\mu\in P(G)$ is called an {\em idempotent measure}
if $\mu*\mu=\mu$. In 1954 Wendel \cite{Wen} proved that each
idempotent measure on a compact topological group coincides with
the Haar measure of some compact subgroup. Later, Wendel's
result was generalized to locally compact groups by Pym \cite{Pym}
and to all topological groups by Tortrat \cite{Tor}. By the {\em Haar
measure} on a compact topological group $G$ we understand the
unique $G$-invariant probability measure on $G$. It is a classical
result that such a measure exists and is unique.  Thus we have the
following characterization of idempotent measures on topological
groups:

\begin{proposition} A probability Radon measure $\mu\in P(G)$ on a topological group $G$
is an idempotent of the semigroup $P(G)$ if and only if $\mu$ is the Haar measure of some compact subgroup of
$G$.
\end{proposition}

We shall use this proposition to describe regular elements of the
convolution semigroups. To this end we apply Proposition 4 of \cite{BHr} that describes
regular elements of the hypersemigroups over topological groups:

\begin{proposition}[Banakh-Hryniv]\label{p4a} For a compact subset $K\in\exp(G)$ of a topological group $G$ the following assertions are equivalent:
\begin{enumerate}
\item $K$ is a regular element of the semigroup $\exp(G)$;
\item $K$ is uniquely invertible in $\exp(G)$;
\item $K=Hx$ for some compact subgroup $H$ of $G$ and some $x\in G$.
\end{enumerate}
\end{proposition}

A similar description of regular elements holds for the convolution semigroup:

\begin{proposition}\label{p4} For a measure $\mu\in P(G)$ on a topological group $G$ the following assertions are equivalent:
\begin{enumerate}
\item $\mu$ is a regular element of the semigroup $P(G)$;
\item $\mu$ uniquely invertible in $P(G)$;
\item $\mu=\lambda*x$ for some idempotent measure $\lambda\in P(G)$ and
some element $x\in G$.
\end{enumerate}
\end{proposition}

\begin{proof} Assume that $\mu$ is a regular element of $P(G)$ and $\nu\in P(G)$ is a measure such that
$\mu*\nu*\mu=\mu$. The measure $\mu*\nu$, being an idempotent of
$P(G)$ coincides with the Haar measure $\lambda$ on some compact
subgroup $H$ of $G$. It follows that $\overline{S_\mu\cdot
S_\nu}=S_{\mu*\nu}=S_\lambda=H$ and hence $S_\mu$ and $S_\nu$ are
compact subsets of the group $G$. Since $\supp:P(G)\to 2^G$ is a
semigroup homomorphism, we get $S_\mu*S_\nu*S_\mu=S_\mu$, which
means that $S_\mu$ is a regular element of the semigroup $\exp(G)$
and hence $S_\mu=\tilde Hx$ for some compact subgroup $\tilde H$
and some element $x\in G$ according to Proposition~\ref{p4a}.

We claim that $\tilde H=H$. Indeed, $H\tilde Hx=S_\lambda
S_\mu=S_{\mu*\nu}S_\mu=S_{\mu*\nu*\mu}=S_\mu=\tilde Hx$ implies
that $H\subset\tilde H$. Next, for any point $y\in S_\nu$ we get
$$\tilde Hxy\subset \tilde HxS_\nu=S_\mu
S_\nu=S_\lambda=H\subset\tilde H$$which yields $xy\in\tilde H$ and
finally $H=\tilde H$.

Next, we show that $\mu=\lambda*x$, which is equivalent to
$\lambda=\mu*x^{-1}$. Observe that $S_{\mu*x^{-1}}=S_\mu
x^{-1}=Hxx^{-1}=H$. Now the equality $\mu*x^{-1}=\lambda$ will
follow as soon as we check that the measure $\mu*x^{-1}$ is
$H$-invariant. Take any point $y\in H$ and note that
$$y*\mu*x^{-1}=y*\mu*\nu*\mu*x^{-1}=y*\lambda*\mu*x^{-1}=\lambda*\mu*x^{-1}=\mu*x^{-1},$$
which means that the measure $\mu*x^{-1}$ on $H$ is
left-invariant. Since $H$ possesses a unique left-invariant
probability measure $\lambda$, we conclude that $\mu=\lambda*x$.

Finally, we show that $\mu$ is uniquely invertible in $P(G)$. It
suffices to check that the measure $\nu$ is equal to
$x^{-1}*\lambda$ provided $\nu=\nu*\mu*\nu$. For this just observe
that $S_\nu$ being a unique inverse of $S_\mu$ is equal to
$x^{-1}H$. Then $S_{x*\nu}=xS_\nu=xx^{-1}H$. Finally, noticing that
for every $y\in H$ we get
$$x*\nu*y=x*\nu*\mu*\nu*y=x*\nu*\lambda*y=x*\nu*\lambda=x*\nu,$$
which means that $x*\nu$ is a right invariant measure on $H$.
Since $\lambda$ is the unique right-invariant measure on $H$ we also
get $x*\nu=\lambda$ and hence $\nu=x^{-1}*\lambda$.
\end{proof}

Given a semigroup $S$ we denote the set of regular elements of $S$ by $\Reg(S)$.

\begin{proposition}\label{p5} For any topological group $G$, the support map $$\supp:\Reg(P(G))\to\Reg(\exp(G))$$is a homeomorphism.
\end{proposition}

\begin{proof} The preceding proposition implies that the map $$\supp:\Reg(P(G))\to\Reg(\exp(G))$$ is bijective.
In order to check the continuity of this map, we must prove that for any open set $U\subset G$ the preimages
$$
\begin{aligned}
\supp^{-1}(U^+)=\;&\{\mu\in \Reg(P(G)):\supp(\mu)\subset U\}\mbox{ and }\\
\supp^{-1}(U^-)=\;&\{\mu\in \Reg(P(G)):\supp(\mu)\cap U\ne\emptyset\}
\end{aligned}
$$are open in $P(G)$. The openness of $\supp^{-1}(U^-)$ follows from the observation that
$\supp(\mu)\cap U\ne\emptyset$ if and only if $\mu(U)>0$. To see that $\supp^{-1}(U^+)$ is open,
fix any measure $\mu\in\Reg(P(G))$ with $\supp(\mu)\subset U$. By Proposition~\ref{p4},
$\supp(\mu)=Hx$ for some compact subgroup $H$ of $G$ and some $x\in G$. The compactness of $H$ allows
us to find an open neighborhood $V$ of the neutral element of $G$ such that $HV^2HV^{-2}HV\subset Ux^{-1}$.
 Now consider the open neighborhood $W=\{\nu\in\Reg(P(G)):\nu(HVx)>\frac12\}$ of the measure $\mu$.
 We claim that $W\subset\supp^{-1}(U^+)$. Indeed, given any measure $\nu\in W$ we can apply Proposition~\ref{p4}
 to find an idempotent measure $\lambda$ and  $y\in G$ such that $\nu=\lambda*y$. Then $\frac12<\nu(HVx)=\lambda(HVxy^{-1})$.
 We claim that $S_\lambda\subset HVVH$. Indeed, given an arbitrary point $z\in S_\lambda$ use the $S_\lambda$-invariance
 of $\lambda$ to conclude that $\lambda(zHVxy^{-1})=\lambda(HVxy^{-1})>1/2$, which implies that the
intersection $zHVxy^{-1}\cap HVxy^{-1}$ is non-empty which yields
$z\in HVxy^{-1}(HVxy^{-1})^{-1}=HVVH$. The inequality
$\lambda(HVxy^{-1})>1/2$ implies that $HVxy^{-1}$ intersects
$S_\lambda$ and hence the set $HVVH$. Then $y\in HV^{-2}HHVx$ and
$S_\nu=S_\lambda*y\subset HV^2HHV^{-2}HVx\subset Ux^{-1}x=U$,
which implies that $\nu\in \supp^{-1}(U^+)$. This completes the
proof of the continuity of the map
$\supp:\Reg(P(G))\to\Reg(\exp(G))$.

The proof of the continuity of the inverse map
$$\supp^{-1}:\Reg(\exp(G))\to\Reg(P(G))$$ is even more involved.
Assume that $\supp^{-1}$ is discontinuous at some point
$K_0\in\Reg(\exp(G))$. By Proposition~\ref{p4a}, $K_0$ is a coset of some compact subgroup of $G$. After a suitable shift, we can assume that $K_0$ is a compact subgroup of $G$ and then $\mu_0=\supp^{-1}(K_0)$ is the unique Haar measure on $K_0$.

Since $\supp^{-1}$ is discontinuous at $K_0$, there is a neighborhood $O(\mu_0)\subset P(G)$ of $\mu_0$ such that $\supp^{-1}(O(K_0))\not\subset O(\mu_0)$ for any neighborhood $O(K_0)\subset\Reg(\exp(G))$ of $K_0$ in $\Reg(\exp(G))$.

It it well-known that the topology of $G$ is generated by the left uniform structure, which is generated by bounded left-invariant pseudometrics. Each bounded left-invariant pseudometric $\rho$ on $G$ induces a
pseudometric $\hat \rho$ on $P(G)$ defined by
$$\hat\rho(\mu_1,\mu_2)=\inf\{\mu(\rho):\mu\in P(G\times G)\;\;P\pr_1(\mu)=\mu_1,\;P\pr_2(\mu)=\mu_2\}$$where $P\pr_i:P(G\times G)\to P(G)$ is the map induced by the projection $\pr_i:G\times G\to G$ onto the $i$th coordinate. By \cite[\S4]{Ba2} or \cite[3.10]{Fe}, the topology of the space $P(G)$ is generated by the pseudometrics $\hat\rho$ where $\rho$ runs over all bounded left-invariant continuous pseudometrics on $G$.

Consequently, we can find a left-invariant continuous pseudometric $\rho$ on $G$ such that the neighborhood $O(\mu_0)$ contains the $\e_0$-ball $B(\mu_0,\e_0)=\{\mu\in P(G):\hat\rho(\mu,\mu_0)<\e_0\}$ for some $\e_0>0$. Replacing $\rho$ by a larger left-invariant pseudometric, we can additionally assume that for the pseudometric space $G_\rho=(G,\rho)$ the map $\gamma: G_\rho\times G_\rho\to G_\rho$, $\gamma:(x,y)\mapsto xy^{-1}$,  is continuous at each point $(x,y)\in K_0\times K_0$ (this follows from the fact that for each continuous left-invariant pseudometric $\rho_1$ on $G$ we can find a continuous left-invariant pseudometric $\rho_2$ on $G$ such that the map $\gamma:G_{\rho_2}\times G_{\rho_2}\to G_{\rho_1}$ is continuous at points of the compact subset $K_0\times K_0$).

The continuity and the left-invariance of the pseudometric $\rho$ implies that the set $G_0=\{x\in G:\rho(x,1)=0\}$ is a closed subgroup of $G$. Let $G'=\{xG_0:x\in G\}$ be the left coset space of $G$ by $G_0$ and $q:G\to G'$, $q:x\mapsto xG_0$, be the quotient projection. The space $G'=G/G_0$ will be considered as a $G$-space endowed with the natural left action of the group $G$. The pseudometric $\rho$ induces a continuous left-invariant metric $\rho'$ on $G'$ such that $\rho(x,y)=\rho'(q(x),q(y))$ for all $x,y\in G$. So, $q:(G,\rho)\to(G',\rho')$ is an isometry.
The pseudometrics $\rho$ and $\rho'$ induce the Hausdorff pseudometrics $\rho_H$ and $\rho_H'$ on the hyperspaces $\exp(G)$ and $\exp(G')$ such that the map $\exp q:\exp(G)\to\exp(G')$ is an isometry. Also these pseudometrics induce the pseudometrics $\hat\rho$ and $\hat\rho'$ on the spaces of measures $P(G)$, $P(G')$ such that the map $Pq: (P(G),\hat\rho)\to (P(G'),\hat\rho')$ is an isometry. The continuity of the map $\gamma:G_\rho^2\to G_\rho$ at $K_0^2$ implies that $(K_0,\rho)$ is a (not necessarily separated) topological group, $K_0\cap G_0$ is a closed normal subgroup of $K_0$ and hence $K_0'=q(K_0)=K_0/K_0\cap G_0$ has the structure of topological group. Then $\mu'_0=Pq(\mu_0)$ is a Haar measure in $K_0'$.

By the choice of the neighborhood $O(\mu_0)$, for every $n\in\IN$ we can find a compact set $K_n\in \Reg(\exp(G))$ such that the measure $\mu_n=\supp^{-1}(K_n)$ does not belong to $O(\mu_0)$. Then $\hat\rho(\mu_n,\mu_0)\ge\e_0$ by the choice of the pseudometric $\rho$.

For every $n\in\IN$ let $\mu_n'=Pq(\mu_n)\in P(G')$,   and  $K'_n=q(K_n)\in\exp(G')$.
The convergence of the sequence $(K_n)$ to $K_0$ in the pseudometric space $(\exp(G),\rho_H)$ implies the convergence of the sequence $(K'_n)$ to
$K'_0$ in the metric space $(\exp(G'),\rho'_H)$, which implies that the union $K'=\bigcup_{n\in\w}K'_n$ is compact in the metric space $(G',\rho')$. Then the subspace $P(K')$ is compact in the metric space $(P(G),\hat\rho')$ and hence the sequence $(\mu_n')_{n\in\IN}$ contains a subsequence that converges to some measure $\mu'$ in $(P(G'),\hat\rho')$. We lose no generality assuming that whole sequence $(\mu_n')_{n\in\IN}$ converges to $\mu'$. Since $\e_0\le\hat\rho(\mu_n,\mu_0)=\hat\rho'(\mu_n',\mu_0')$, we conclude that $\mu'\ne\mu_0'$. We shall derive a contradiction (with the uniqueness of a left-invariant probability measure on compact groups) by  showing that $\mu'$ is a left-invariant measure on $K_0'$, distinct from the Haar measure $\mu_0'$.

The $\hat\rho'$-convergence $\mu_n'\to\mu'$ and $\rho'_H$-convergence $\supp(\mu'_n)=K_n'\to K_0'$ imply that $\supp(\mu')\subset K_0'$ and thus $\mu'$ is a probability measure on the compact topological group $K_0'$. It remains to check that the measure $\mu'$ is left-invariant. Assuming the converse, we can find a point $a\in K_0'$ such that $a*\mu'\ne\mu'$ and thus $\e=\hat\rho'(\mu',a*\mu')>0$. Since the map $\gamma:G_\rho\times G_\rho\to G_\rho$ is continuous at each point $(x,y)\in K_0\times K_0$, we can find
a positive $\delta<\frac\e4$ so small that $\rho(xy,x'y)<\frac\e4$ for any $x,y\in K_0$ and $x'\in G$ with $\rho(x',x)<\delta$. Since $\rho_H(K_n,K_0)\to 0$ and $\hat\rho'(\mu_n',\mu')\to0$, there is a number $n\in\IN$ and a point $a_n\in K_n$ such that $\rho(a,a_n)<\delta$ and $\hat\rho'(\mu'_n,\mu')\le\e/4$.
Consider two left shifts $l_a:G\to G$, $l_a:x\mapsto ax$, and $l_{a_n}:G\to G$. The choice of $\delta$ guarantees that $\rho_{K_0}(l_a,l_{a_n})=\sup_{x\in K_0}\rho(l_a(x),l_{a_n}(x))\le \frac\e4$. Then $$\hat\rho'(a*\mu',a_n*\mu')=\hat\rho'(Pl_a(\mu'),Pl_{a_n}(\mu'))\le\frac\e4.$$
The left shift $l_{a_n}:G\to G$, being an isometry of the pseudometric space $(G,\rho)$, induces an isometry $l'_{a_n}:G'\to G'$ of the metric space $(G',\rho')$, which induces the isometry $Pl'_{a_n}:P(G')\to P(G')$ of the corresponding space of measures. So, $\hat\rho'(a_n*\mu',a_n*\mu'_n)=\hat\rho'(Pl'_{a_n}(\mu'),Pl'_{a_n}(\mu'_n))=\hat\rho'(\mu',\mu'_n)\le\frac\e4$. The compact set $K_n$, being a regular element of the semigroup $\exp(G)$ is equal to $H_nx_n$ for some compact subgroup $H_n\subset G$ and some point $x_n\in G$ according to Proposition~\ref{p4a}.
Then $\mu_n=\supp^{-1}(K_n)$ is equal to $\lambda_n*x_n$ where $\lambda_n$ is the Haar measure on the group $H_n$. Since $\lambda_n$ is left-invariant, $a_n*\mu_n=a_n*\lambda_n*x_n=\lambda_n*x_n=\mu_n$ and hence $a_n*\mu_n'=\mu'_n$.

Now we see that
$$
\begin{aligned}
\hat\rho'(\mu',a*\mu')&\le \hat\rho'(\mu',\mu'_n)+\hat\rho'(\mu'_n,a_n{*}\mu'_n)+
\hat\rho'(a_n{*}\mu'_n,a_n{*}\mu')+\hat\rho'(a_n{*}\mu',a{*}\mu')\le\\
&\le \frac{\e}4+0+\frac{\e}4+\frac{\e}4<\e=\hat\rho'(\mu',a*\mu'),
\end{aligned}$$which is a desired contradiction.
\end{proof}

The following corollary establishes the first part of
Proposition~\ref{pr1}. The second part of that proposition follows
from Theorem 2 of \cite{BHr}.

\begin{corollary} Let $G$ be a topological group. Then a topological regular semigroup $S$ can be embedded into the
 hypersemigroup $\exp(G)$ if and only if $S$ can be embedded into the convolution semigroup
$P(G)$.
\end{corollary}

\begin{proof}  If $S\subset\exp(G)$ is a regular subsemigroup, then $S\subset \Reg(\exp(G))$ and $\supp^{-1}(S)$
is an isomorphic copy of $S$ in $P(G)$ according to Propositions~\ref{p5}.  Conversely, if $S\subset P(G)$ is a regular
subsemigroup, then its image $\supp(S)$ is an isomorphic copy of $S$ in $\exp(G)$.
\end{proof}

\section*{Acknowledgements}
This research was supported by the
Slovenian Research Agency grants P1-0292-0101-04, J1-9643-0101 and
J1-2057-0101. We thank the referee for comments and suggestions.

\end{document}